\documentclass{amsart}
\usepackage{geometry}                
\usepackage{graphicx}
\usepackage{amssymb}
\usepackage{epstopdf}
\newtheorem{theorem}{Theorem}
\newtheorem{prop}[theorem]{Proposition}

\newtheorem{lemma}[theorem]{Lemma}
\newtheorem{cor}[theorem]{Corollary}
\theoremstyle{definition}
\newtheorem{defin}[theorem]{Definition}
\newtheorem{remark}[theorem]{Remark}

\newcommand{\idiot}[1]{\vspace{5 mm}\par \noindent
\marginpar{\textsc{Note}}
\framebox{\begin{minipage}[c]{0.95 \textwidth}
#1 \end{minipage}}\vspace{5 mm}\par}

\renewcommand{\idiot}[1]{}

\def\myem{\em}
\def\Gasc{{C}}
\def\GAsc{{\Gasc}}
\def\Sn{{\mathfrak S}}

\def\B{{\mathfrak B}}
\def\D{{\mathfrak D}}

\def\SS{{S}}
\def\PP{{\mathcal P}}
\def\GG{{G}}
\def\FF{{F}}
\def\std{{\textnormal{std}}} 
\def\Des{{\textnormal{Des}}}
\def\GDes{{\textnormal{GDes}}}

\title{Posets related to the connectivity set of Coxeter groups}

\author{N. Bergeron}\address[Nantel Bergeron]{Department of Mathematics and Statistics\\ York  University\\
   To\-ron\-to, Ontario M3J 1P3\\ CANADA} \email{bergeron@mathstat.yorku.ca}  \urladdr{http://www.math.yorku.ca/bergeron}

\author{C. Hohlweg}\address[Christophe Hohlweg]{The Fields Institute\\ 222 College Street\\
Toronto, Ontario, M5T 3J1\\ CANADA} \email{chohlweg@fields.utoronto.ca}
\urladdr{http://www.fields.utoronto.ca/\~{}chohlweg}

\author{M. Zabrocki} \address[Mike Zabrocki]{Department of Mathematics and Statistics\\ York  University\\
    To\-ron\-to, Ontario M3J 1P3\\ CANADA} \email{zabrocki@mathstat.yorku.ca}  \urladdr{http://www.math.yorku.ca/\~{}zabrocki}

\date{\today}

\thanks{This work is supported in part by CRC and NSERC.}

\begin{document}
\maketitle

\begin{abstract}
We define the notion of connectivity set for elements of any finitely
generated Coxeter group. 
Then we define an order related to this new statistic and show that
the poset is graded and each interval is a shellable lattice. This
implies that any interval is Cohen-Macauley. We also give a Galois
connection between intervals in this poset and a boolean poset. This
allows us to compute the M\"obius function for any interval.
\end{abstract}

%
%
%
%

\section*{Introduction}\label{se:Intro}

Let $\Sn_n$ denote the symmetric group, the group of  permutations on the set $[n] = \{1,2,\ldots, n\}$.
 The length $\ell(w)$ of a permutation $w\in \Sn_n$ is its number of inversions i.e. the number
 of pair $(i,j)$ with $1\leq i<j\leq n$ and $w(i)>w(j)$.
The {\em connectivity set}  (also known as {\em global ascent set}) of a permutation $w = w(1) w(2) \cdots w(n) \in \Sn_n$ is equal to
\begin{equation}\label{gasc:typeA}
\Gasc(w) = \{ i\in [n-1] : w(j) < w(k), \forall 1\leq j \leq i < k\leq n \}.
\end{equation}
This set has recently been the subject of an article by R.~Stanley  \cite{St} which
is related to the number of connected components of a permutation (see \cite{calan,C}).
It turns out that the connectivity set appears also to be linked with the
study of two combinatorial Hopf algebras:
the Malvenuto-Reutenauer Hopf algebra and  the quasisymmetric
functions in non-commuting variables \cite{AS,BZ,NCSF6}.

In  their analysis of the Malvenuto-Reutenauer Hopf algebra of permutations \cite{AS},
 Aguiar and Sottile   consider the set $\GDes(w)$ of global descents of a permutation $w\in \Sn_n$
 which is related to the connectivity set by the formula
 $\GDes(w) = w_0 (\Gasc(ww_0))$ where $w_0=n\,n-1\,\cdots\,2\,1$ is the unique permutation
 of maximal length in $\Sn_n$. Similar notions are also used in \cite{NCSF6}.  

\begin{figure}[htbp] 
   \centering
   \includegraphics[width=2.5in]{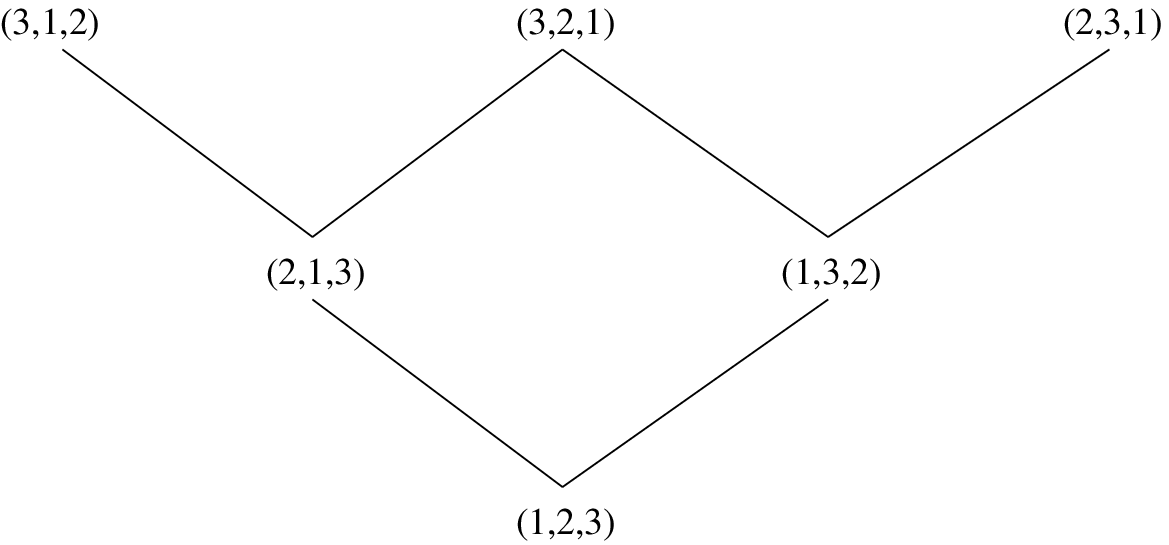}
   \caption{Hasse diagram for posets of $(\Sn_3, \leq)$}
   \label{hasseS3}
\end{figure}

\begin{figure}[htbp] 
   \centering
   \includegraphics[width=5.5in]{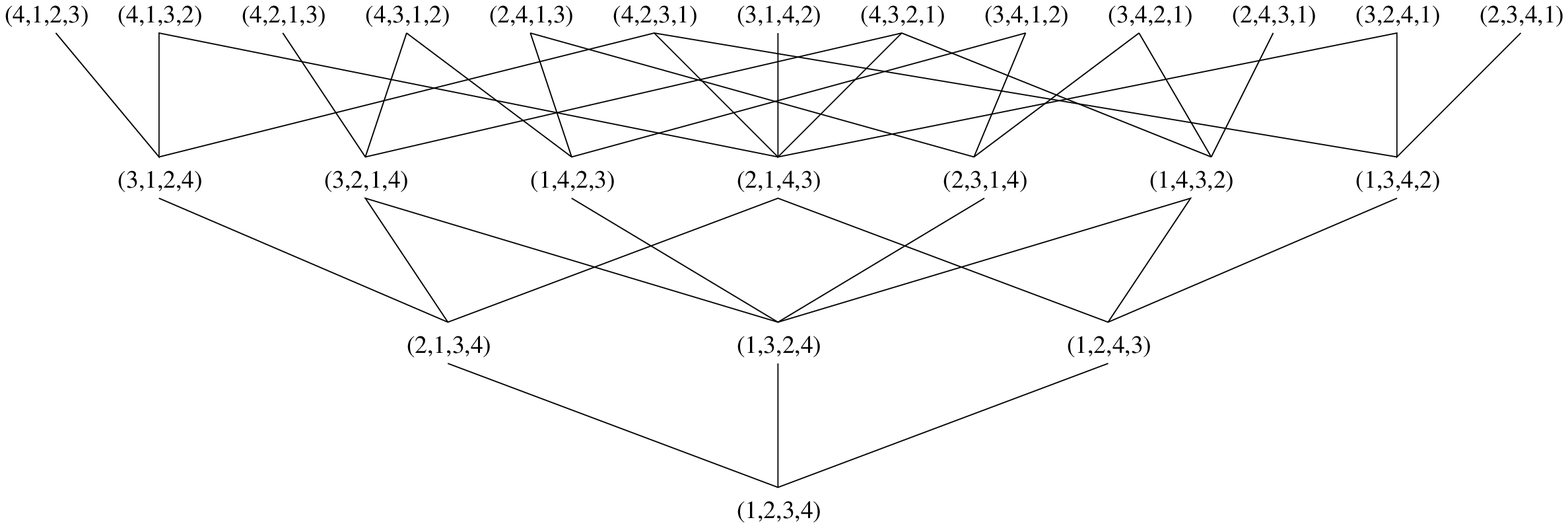}
   \caption{Hasse diagram for posets of $(\Sn_4, \leq)$}
   \label{hasseS4}
\end{figure}

In the  work of Bergeron and Zabrocki \cite{BZ} a natural order on set compositions
arose out of the Hopf algebra structure of the quasisymmetric functions in
non-commuting variables.  This order came out of some simple conditions
placed on the comultiplicative structure of the Hopf algebra.  Restricting
attention to set compositions with one element in each part, the partial order
can be viewed as a partial order $\leq$ on permutations: let $u,v\in\Sn_n$, then 
\begin{equation}\label{order:typeA}
u\leq v\hbox{ if and only
 if }\Gasc(v)\subset \Gasc (u)\hbox{ and }\std_{\Gasc(u)}(v)=\std_{\Gasc(u)}(u).
 \end{equation} Here $\std_{I}(w)$
 denotes the standardization of $w\in\Sn_n$ along the blocks of $I\subset [n-1]$. For instance
$\std_{\{1,4\}}(4.623.51)=1.423.65$. 
A few of the Hasse diagrams for these posets appear in
Figures~\ref{hasseS3} and~\ref{hasseS4}.

The study of the poset $(\Sn_n,\leq)$ led us to consider a
generalization of these concepts in the language of Coxeter
groups. Not only is the definition of this poset strongly
simplified,
but the connectivity set turns out 
to be closely linked with a classical result in the theory of
Coxeter groups: {\myem the word property}, due to Tits \cite{tits} and independently to Matsumoto
\cite{matsumoto} (see for instance \cite[Theorem 3.3.1]{bjorner-brenti}).\\

Let $(W,S)$ be a finitely generated Coxeter system whose length function is denoted by $\ell:W\to\mathbb N$.
 Denote by $e$ the identity of $W$.
 We denote by $W_I$ the {\myem parabolic subgroup} of $W$ generated by $I\subset S$.
Let $w\in W$,
 the word property says that any pair of reduced expressions for $w$  can be linked by a sequence
of braid relation transformations. In particular, the set
\begin{equation}
S(w)=\{s_i\in S\,:\, w= s_1\dots s_{\ell(w)} \hbox{ reduced}\}
=\bigcap_{{\scriptstyle I\subset S \atop \scriptstyle w\in W_I}}I
\end{equation}
is independent of the choice of a reduced expression for $w$. It is clear that $w\in W_{S(w)}$.
 
The {\myem descent set} of $w\in W$ is the set
$$
\Des(w)=\{s\in S : \ell(ws)<\ell(w)\}.
$$
Let $I\subset S$,  it is well-known that the set
$$
  X_I =\{u\in W : \ell(us)>\ell(u),\,\forall s\in I\} =\{ u \in W : \Des(u) \subseteq S\backslash I \}
$$
is a set of minimal length coset representatives of $W/W_I$.  Each element  $w \in W$ has a
unique decomposition $w = w^I w_I$ where $w^I \in X_I$ and $w_I \in W_I$ and moreover 
$\ell(w) = \ell(w^I) + \ell(w_I)$.  The pair $(w^I, w_I)$ is generally referred to as the
{\myem parabolic components of $w$ along $I$} 
(see \cite[Proposition 2.4.4]{bjorner-brenti} or \cite[5.12]{H}).

\begin{defin} $\ $
\begin{enumerate}
\item  The {\myem connectivity set of $w\in W$} is the set  $\Gasc(w)=S\setminus S(w)$.
\item Let $u,v\in W$, $u\leq v$ if and only if the parabolic component $v_{S(u)}=u$.
\end{enumerate}
\end{defin}

In other words, $\Gasc(w)$ is the set  of simple reflections which do not 
appear in a reduced expression for $w$.
We will show in \S\ref{se:symgrp} that in the case where $W=\Sn_n$ and $S$ is the set  of simple transpositions
$\tau_i = (i, i+1)$, with $i\in [n-1]$, the definitions coincide with those in equations 
\eqref{gasc:typeA} and \eqref{order:typeA}.\\

Let $P$ be a poset and $u\in P$, a {\myem cover} of $U$ is an element $v>u$ in $P$ such
 that the interval $[u,v]=\{u,v\}$. 
Recall that a lattice $L$ is called {\myem upper semimodular} if for
$w,g,h\in L$ such that $g$ and $h$ cover $w$, there is a $x$
which covers both $g$ and $h$. 

Our main results is the following theorem.

\begin{theorem}\label{thm:Main} Let $(W,S)$ be a finitely generated Coxeter system.
\begin{enumerate}
\item  The poset $(W,\leq)$ is graded. The rank function is $w\mapsto |S(w)|$.
\item The interval $[u,v]$ is an upper semimodular lattice, for any
$u\leq v$ in $W$.
\end{enumerate}
\end{theorem}

We obtain immediately the following corollary.  We refer the reader 
to \cite{Bj,St1,St2} and \cite[Appendix A2]{bjorner-brenti} for more 
information about the concepts of shellability and Cohen-Macauliness.

\begin{cor}\label{cor:Main}
The order complex of $[u,v]$  is  shellable (hence
Cohen-Macaulay), for any $u\leq v$ in $W$.
\end{cor}

 In \S\ref{se:Main}, we construct a {\myem Galois connection}
 between any interval $[u,v]$ and a boolean
poset. This will allow us to show our second main result.

\begin{theorem} \label{thm:mobius} For any $u\le v$ in $W$, the M\"obius function is given by
$$\mu(u,v) = \begin{cases}(-1)^{|\SS(v)|-|\SS(u)|}&\hbox{ if }\SS(u) = S(v) \backslash \Des(v^{S(u)})
\hbox{ or }u=v\\
0&\hbox{ otherwise.}
\end{cases}$$
In particular,
$$\mu(e, v) = \begin{cases} (-1)^{|\SS(v)|} &\hbox{if }v = w_0(\SS(v))\\
0&\hbox{otherwise,}
\end{cases}$$
where $w_0(I)$ denote the element of maximal length in $W_I$.
\end{theorem}

 It is interesting to remark that this M\"obius function in type $A$ correspond to the coefficients of the primitives of the Malvenuto-Reutenauer Hopf algebra. This was done by  Duchamp, Hivert and Thibon \cite[\S 3.4]{NCSF6}. They defined operations $p_I(w)$  for a permutation $w$ that correspond to our $\std_I(w)$, but they do not consider the order it induces.

Finally in \S\ref{se:Examples}, we will  give
a characterization of the connectivity set when the Coxeter group is finite  of type $B$ and $D$
 and when they are considered as subgroups of the symmetric group acting on the set $(-[n])\cup[n]$.
  We also give formulas for generating functions for the numbers of
elements of a Coxeter group of type $A$, $B$ and $D$ with exactly $k$ elements in the connectivity set.

\bigskip\bigskip
\section*{Acknowledgment}
\begin{minipage}[c]{0.95 \textwidth}
The authors would like to thank Vic Reiner for helpful
suggestions on this research. 
He has outlined a proof that the poset $W - \{e\}$
is contractible and remarked that
the order complex of $(W,\le)$ was shellable for
$W = \Sn_n$ and small values of $n$. It is still open to
show this in general, since $(W,\le)$ is not a single interval.
We also thank Hugh Thomas for valuable discussions.
\end{minipage}

\idiot{{\tt 
On 14-Sep-05, at 6:38 PM, Victor Reiner wrote:

Hi Mike,

  I do have two comments about your paper and what you asked.

(1) I think I see now that there is an easy argument for the contractibility of the poset on $W - \{e\}$.

You first need to know the following commonly used special case of Quillen's Fiber Lemma:

 if the subposet $P_{<x}$ consisting of all elements strictly below x
 in some poset P is contractible, then P is homotopy equivalent to $P-\{x\}$.

Secondly, you need to prove a strengthening of your Theorem 4, characterizing the homotopy of an arbitrary (open) interval (x,y) as either spherical or contractible.  I think this already follows from your Galois connection;  I didn't look at it too closely.

Lastly, you notice that elements x in the interval $[e,w_0]$ are
exactly the ones that have S(x) = Des(x), and hence they are the ones
for which $(e,x)$=$P_{<x}$ is contractible, so that they can be
removed by the above lemma, without affecting the homotopy type
(and you remove all such x one by one starting from the higher ones
first, then the lower ones).  At the end you're left with the semiopen
interval $(e, w_0]$, which has $w_0$ as a top element, so it's contractible.

(2) The following enumerative problem looks tractable to me (and I remember thinking about doing it long ago, but I didn't know if anyone cared about the connectivity sets at that time):  go back to my old paper "The distribution of descents and length in a Coxeter group", and instead of just computing recursively the joint distribution of

  ( length l(w),  descent set Des(w) )

as I did there, do the more general calculation of

  ( length l(w),  descent set Des(w), connectivity set C(w) )

or perhaps replace connectivity set by its complement $S\backslash C(w)$. I think this is do-able.

Best wishes,
Vic}
}
%
%
%
%
%

\section{Main Results}\label{se:Main}

For additional information about finitely generated Coxeter groups
we refer the reader to \cite{H,bjorner-brenti}. In this section,
$(W,S)$ is an arbitrary finitely generated Coxeter system (except
at the end of \S\ref{se:symgrp} where we will discuss the case of
symmetric groups).

\subsection{Preliminaries}\label{se:Prelim}

We start by recalling a well-known lemma whose proof follows
 immediately from definitions.

\begin{lemma} \label{lemma:weakGasc}
For $u,v \in W$ such that $\ell(u v) = \ell(u) + \ell(v)$, we have
$\Gasc(uv) = \Gasc(u) \cap \Gasc(v)$ (or alternatively, $\SS(uv) =
\SS(u) \cup \SS(v)$).
\end{lemma}

The {\myem (left) weak order} $\leq_\mathcal{L}$ on $W$ may be defined as
follows: $u\leq_\mathcal{L} v$ if and only if there is $v'\in W$
such that $v=v'u$ and $\ell(v)=\ell(v')+\ell(u)$. 
\begin{prop}\label{prop:EasyConsequences} Let $u,v\in W$.

\begin{enumerate}

\item If $u\leq v$ in $W$, then $\Gasc(v)\subset \Gasc(u)$ (or
alternatively, $S(u)\subset S (v)$).

\item For $u\leq v$ in $W$, $[u,v]=\{v_I : S(u)\subset I\subset
S(v)\}$.

\item if $u\leq v$, then $u \leq_\mathcal{L} v$.

\item If $W_{S(u)}$ is of finite index $k$ in $W$, then
 $|\{w\in W : w\geq u\}|=k$.

\end{enumerate}

\end{prop}

\begin{proof} (1)
 As $v = v^{\SS(u)} u$ with $\ell(v) = \ell(v^{\SS(u)}) +
\ell(u)$, then  $\Gasc(v) = \Gasc(v^{\SS(u)}) \cap \Gasc(u)
\subseteq \Gasc(u)$ by Lemma~\ref{lemma:weakGasc}.

(2) and (3) follow from definitions.

(4) There is a bijection between the set $X_{\SS(u)}$ and $\{ w
\in W | w \geq u \}$.  $w \geq u$ if and only if $w = w^{\SS(u)}
u$, hence the map
\begin{equation*}
\{ w \in W | w \geq u \} \longrightarrow X_{\SS(u)}
\end{equation*}
which sends $w \mapsto w^{\SS(u)}$ has as inverse $x \mapsto xu$
and so is a bijection.
\end{proof}

\begin{remark} In fact, we could define the order on the 
 set of coset $W/W_I$, for any $I\subset S$. Indeed,  we can take
 an element $w$ such that $C(w)=I$. Then by
Proposition~\ref{prop:EasyConsequences} (4) the poset
 $(W/W_I,\leq)$ is obtained by taking the set
 $\{gW_I : g\geq w\}$ and the induced order. 
For  the symmetric group ($W = \Sn_n$), this gives the order on set compositions
 considered in \cite{BZ}.
\end{remark}

Now we recall some useful facts about minimal coset
representatives. If $I \subseteq J \subseteq S$, then $X_I^J = X_I
\cap W_J$ is the set of minimal coset representatives of $W_J
\slash W_I$ and $X_J X_I^J = X_I$. The following lemma is
well-known.

\begin{lemma}\label{lem:Useful} Let $I\subset J\subset S$ and $w\in W$, then
$w_I=(w_J)_I$. In particular, $w_I \leq w_J$.
\end{lemma}
\begin{proof}
Write $w=w^Iw_I=w^Jw_J$ and $w_J=(w_J)^I(w_J)_I$ and conclude by
uniqueness of the parabolic components.
\end{proof}

Let $K \subseteq S$.  If $W_K$ is finite, then $W_K$ contains a
unique element $w_0(K)$ of maximal length. It can be characterized
 as the unique element $w$ in $W_K$ such that $\Des(w)=K$.
 In fact, for any $w\in W$, $W_{\Des(w)}$ is finite (see for instance
 \cite[Proposition 2.3.1]{bjorner-brenti}).

\begin{lemma}\label{lem:W0} Let $I\subset J\subset S$ such that
$W_J$ is finite, then $(w_0(J))_I=w_0(I)$.
\end{lemma}
\begin{proof} Assume there is an $s\in I$ such that $\ell((w_0(J))_I
s)>\ell((w_0(J))_I)$. As $((w_0(J))^I,(w_0(J))_Is)$ are the
parabolic components of $w_0(J)s$ along $I$ we obtain
$$
\ell(w_0(J)s)=\ell((w_0(J))^I)+\ell((w_0(J))_Is)>\ell(w_0(J)).
$$
Hence $s\notin Des(w_0(J))=J$ which contradicts $s\in I\subset J$.
 Therefore $Des((w_0(J))_I)=I$.
\end{proof}

\smallskip

\subsection{The symmetric group} \label{se:symgrp}
We end this preliminary discussion by proving the equivalence of
definitions on the case of $W = \Sn_n$ and $S$ is the set  of
simple transpositions $\tau_i = (i, i+1)$, with $i\in [n-1]$.

The {\sl standardization} of a word $w=a_1 a_2\cdots a_n$ of
length $n$ in an totally ordered alphabet, denoted
 by $\std(w)$, is the unique permutation $\sigma\in\Sn_n$
 such that for all $i<j$
 we have $\sigma(i)>\sigma(j)$ if and only if $a_i>a_j$.
For instance, for $w=cbbcaa$ in the alphabet $\{a<b<c\}$ we have
$\std(w)=534612$. A {\sl composition} of $n$ is a sequence
$\mathbf c=(c_1,\dots,c_k)$ of positive integers whose sum is $n$.
There is a well-known bijection between compositions of $n$ and
subsets of $[n-1]$ defined by
\begin{eqnarray*}
I=\{i_1,i_2,\dots,i_k\} &\mapsto& \mathbf c_I
=(i_1,i_2-i_1,\dots,i_k-i_{k-1},n-i_{k}).
\end{eqnarray*}

 Let $I\subset S$ and $\mathbf c_{S\setminus I}=(c_1,c_2,\ldots,c_k)$.
  Set $t_i=c_1+c_2+\cdots+c_i$ for all $i$. Given a $k$-tuple
$(\sigma_1, \sigma_2,\ldots,\sigma_k)\in\Sn_{c_1}\times\Sn_{c_2}
\times\cdots\times\Sn_{c_k}$ of permutations, we define
$\sigma_1\times\sigma_2\times\cdots\times\sigma_k\in\Sn_n$ as the
permutation that maps an element $a$ belonging to the interval
$[t_{i-1}+1,t_i]$ onto $t_{i-1}+\sigma_i(a-t_{i-1})$. This
assignment defines an isomorphism $\mathfrak
S_{c_1}\times\mathfrak S_{c_2}\times\cdots \times\mathfrak
S_{c_k}\simeq W_{I}$. We have also the well-known characterization
$$
  X_{I}=\bigl\{\sigma\in\mathfrak S_n\bigm|\forall i,\ \sigma
\text{ is increasing on the interval }[t_{i-1}+1,t_i]\bigr\}.
$$
Specifically, write $\sigma\in\Sn_n$ as
 the concatenation  $\sigma_1\cdot\ldots\cdot\sigma_k$ of words in the alphabet $\mathbb N$
 such that  the length of the word $\sigma_i$ is $c_i$. It is then easy to check that
  the parabolic components of $\sigma$ along $I$ are
\begin{equation}
\label{eq:ParabolicCompenent}%
 \sigma_{I}=\std_{S\setminus I}(\sigma):=\std(\sigma_1)\times
\std(\sigma_2)\times\cdots\times\std(\sigma_k) \in \Sn_{\mathbf c}
\quad\text{and}\quad \sigma^{I}=\sigma \sigma_{S \setminus
I}^{-1}\in X_{I}.
\end{equation}

\begin{prop}  For $u \in \Sn_n$,
$$
\Gasc (u) = \{ \tau_i \in S : u(j) < u(k), \forall~1\leq j \leq i
< k<n \}.
$$
\end{prop}
\begin{proof}
For $u \in \Sn_n$, write $\Gasc ( u) = \{ \tau_{i_1}, \tau_{i_2},
\ldots, \tau_{i_k} \}$ with $1 \leq i_1 < i_2 < \cdots < i_k \leq
n-1$.  Then we have that
$$u \in
W_{S(u)} \simeq \mathfrak
S_{c_1}\times\mathfrak S_{c_2}\times\cdots \times\mathfrak
S_{c_k}.
$$
If $\tau_i \in \GAsc (u)$, then for each $j < i$, $u(j) \leq i$
and for each $k>i$, $u(k) > i$.  Hence $\GAsc(u) \subseteq \{
\tau_i : u(j) < u(k), \forall j \leq i < k \}$.  Conversely, if
$\tau_k \in  \{ \tau_i : u(j) < u(k), \forall j \leq i < k \}$,
then $u \in \Sn_{k} \times \Sn_{n-k}$. Hence $\tau_k \notin S(u)$
implying $\tau_k\in \Gasc(u)$.
\end{proof}

As $u_{S(u)}=u$ and by (\ref{eq:ParabolicCompenent}) and
Proposition~\ref{prop:EasyConsequences} (1), we obtain immediately
the following corollary.

\begin{cor} Let $u,v\in\Sn_n$, the following propositions are
equivalent:
\begin{enumerate}
\item $u\leq v$;
\item $\std_{\Gasc(u)}(v)=u$;
\item $\Gasc(v)\subset \Gasc(u)$ and
$\std_{\Gasc(u)}(v)=\std_{\Gasc(u)}(u)$.
\end{enumerate}
\end{cor}

\subsection{Proof of Theorem~\ref{thm:Main}}\label{sse:ProofMain}

From the examples given in the previous section, we see
that this poset does not form a lattice, but each interval
$[u, v]$ in this poset does.

\begin{lemma} \label{lem:latticedef} Let $u,v \in W$ such that $u < v$.  For any $w,g \in [u,v]$
we have
\begin{enumerate}
\item[$(a)$] $v_{\SS(w) \cup \SS(g)}$ is the unique
least upper bound of $g$ and $w$
in $[u,v]$ (i.e. $w \vee g := \min_\leq\{ x : v \geq x \geq w, v
\geq x \geq g \} = v_{\SS(w) \cup \SS(g)}$). Moreover, $\SS(w \vee
g) = \SS(w) \cup \SS(g)$. \item[$(b)$] $v_{\SS(w) \cap \SS(g)}$ is
the unique greatest lower bound of $w$ and $g$ in $[u,v]$ (i.e. $w \wedge g := \max_\leq\{ x: u \leq x
\leq w, u \leq x \leq g \} = v_{\SS(w) \cap \SS(g)}$).  Moreover,
$\SS(w \wedge g) \subseteq \SS(w) \cap \SS(g)$.
\end{enumerate}
\end{lemma}

\begin{proof}
$\SS(w) \subseteq \SS(w) \cup \SS(g)$, so by
Lemma~\ref{lem:Useful}, $w = v_{\SS(w)} \leq v_{\SS(w) \cup
\SS(g)}$ and hence $\SS(w) \subseteq \SS(v_{\SS(w) \cup \SS(g)})$.
Similarly, $\SS(g) \subseteq \SS(v_{\SS(w) \cup \SS(g)})$.
From this, we have $\SS(w) \cup \SS(g) \subseteq \SS(v_{\SS(w) \cup
\SS(g)}) \subseteq \SS(w) \cup \SS(g)$ and therefore $\SS(v_{\SS(w) \cup
\SS(g)}) =\SS(w) \cup \SS(g)$. 
Conversly, if $h \geq w, g$, then
$\SS(h) \supseteq \SS(w) \cup \SS(g)$ and so $h = v_{\SS(h)} \geq
v_{\SS(w) \cup \SS(g)}$.

For part $(b)$, observe that $\SS( v_{\SS(w) \cap \SS(g)} )
\subseteq \SS(w) \cap \SS(g)$ since $v_{\SS(w) \cap \SS(g)} \in
W_{\SS(w) \cap \SS(g)}$. Now $\SS(w) \cap \SS(g) \subseteq
\SS(w), \SS(g)$ gives $v_{\SS(w) \cap \SS(g)} \leq v_{\SS(w)} =w,
v_{\SS(g)} = g$.  Moreover, if $h \leq g,w$, then $\SS(h)
\subseteq \SS(w) \cap \SS(g)$ and hence $h = v_{\SS(h)} \leq
v_{\SS(w) \cap \SS(g)}$.
\end{proof}


\begin{remark}
We do not have $\SS(w \wedge g) = \SS(w) \cap \SS(g)$ in general.
For instance, in $\Sn_4$,  the two elements $s_2 s_1$ and $s_2 s_3$.  $s_2 s_1
\wedge s_2 s_3 = e$ and $\SS(s_2 s_1) \cap \SS(s_2 s_3) = \{ s_2
\}$.
\end{remark}


\begin{lemma}\label{lem:findintermediates}
Let $u,v \in W$ be such that $u < v$.  For any $s \in
\Des(v^{\SS(u)})$ we have
$$u < (v_{\SS(u) \cup \{s\}})^{\SS(u)} u = v_{\SS(u) \cup \{ s\}} < v.$$
Moreover, $\Gasc((v_{\SS(u) \cup \{s\}})^{\SS(u)} u) = \Gasc(u)
\backslash \{s\}$.
\end{lemma}

\begin{proof}  Let $x = (v_{\SS(u) \cup \{ s\}})^{\SS(u)}$ so that
$u = v_{\SS(u)}  \leq v_{\SS(u) \cup \{ s\}} = xu$.  Choose $s \in
\Des( v^{\SS(u)} )$, then $s \notin \SS(u)$ since $v^{\SS(u)} \in
X_{\SS(u)} = \{ w \in W : \Des(w) \subseteq S \backslash \SS(w)
\}$. As $x\in X_{S(u)}$, any reduced expression of $x$ must
end in a simple reflection $r\in S\setminus S(u)$ (see for
instance \cite[Lemma 2.4.3]{bjorner-brenti}). But $x \in
W_{S(u)\cup \{s\}}$, therefore $r\in S(u)\cup \{s\}\setminus S(u)$.
In other words, any reduced expression of $x$ must end in $s$.
Hence $s\in S(x)$.
 By Lemma \ref{lemma:weakGasc}, $\SS(xu) = \SS(u) \cup \SS(x)
= \SS(u) \cup \{ s \}$.

This implies that $\Gasc(xu) = \Gasc(u) \backslash \{ s\}$.
Moreover, $v = v^{\SS(u) \cup \{ s \}} xu$ and $xu \in W_{\SS(u)
\cup \{ s \}}$. By the uniqueness of the parabolic components, we
have that $v_{\SS(u) \cup \{ s \}} = xu$ and so can conclude that
$xu \leq v$.
\end{proof}

The next proposition identifies exactly which elements cover $u$.

\begin{prop}\label{prop:Covers}
Let $u < v$.  If $u'$ is a cover of $u$ in $[u,v]$, then $u' =
(v^{\SS(u)})_{\SS(u) \cup \{ s \}} u$ for some $s \in
\Des(v^{\SS(u)})$.  Moreover, $\Gasc(u') = \Gasc(u) \backslash \{
s \}$.
\end{prop}

\begin{proof} Let $u<u'$ be a cover.
By Lemma \ref{lem:findintermediates}, there is an $s \in
\Des((u')^{\SS(u)})$ and an $x = ( (u')_{\SS(u) \cup \{ s \}}
)^{\SS(u)}$ such that $x \neq id$ and $u < xu < u'$. Since $u'$ is
a cover of $u$, then $xu = u'$.  Since $v = v^{\SS(u) \cup \{ s\}}
x u$ with $\ell(v) = \ell(v^{\SS(u) \cup \{ s\}}) + \ell(x) +
\ell(u)$, therefore $v^{S(u)} = v^{\SS(u) \cup \{ s\}} x$ and $s
\in \Des(v^{S(u)})$. Hence $x=(v^{\SS(u)})_{\SS(u) \cup \{ s \}} $ as expected.
\end{proof}

\begin{cor}
$$| \{ u' : u'\hbox{ covers }u\hbox{ in }[u,v] \}| = |\Des(v^{\SS(u)})|.$$
\end{cor}

\begin{remark}
The poset $(W, \leq)$ does not not form a lattice, even after
adding a maximal element, since $u \vee v$ does not exist in
general (see for example the Hasse diagram for $( \Sn_4, \leq)$,
where the two elements $(1,4,2,3)$ and $(2,3,1,4)$ do not have a
unique least upper bound). In addition, the intervals are not in general
distributive lattices (see for example the interval $[(1,2,3,4),
(4,2,3,1)]$ in the poset $( \Sn_4, \leq)$).
Moreover an interval is not modular since is not lower semimodular.
 Take in $\Sn_4$ the interval $[1234,4231]$. Then $4231$ covers
 both $3124$ and $1342$ but $3124\wedge 1342=1234$.
\end{remark}

\begin{proof}[Proof of Theorem~\ref{thm:Main}] Choose $u,v \in W$ such that $u<v$.

(1) By Lemma \ref{lem:findintermediates} we can construct a path
$u_0=u < u_1 < u_2 < \cdots < u_k < v = u_{k+1}$ such that
$|\SS(u_{i+1})| = |\SS(u_{i})| +1$.  Since if $u'$ covers $u$,
then it also holds that $|\SS(u')| = |\SS(u)|+1$, hence the poset
is graded.

(2) The fact that $([u,v],\wedge,\vee)$ is a lattice follows from
Lemma~\ref{lem:latticedef}.
Now take $w\in [u,v]$ and $g,h\in [u,v]$ covering $w$. By
Proposition~\ref{prop:Covers} we know that there are distinct
$s,r\in \Des(v^{S(w)})$ such that $S(g)=S(w)\cup \{s\}$ and
$S(h)=S(w)\cup \{r\}$. By Lemma~\ref{lem:latticedef} we know that
 $S(g\vee h)=S(w)\cup \{r,s\}$. In other words $g\vee h$ covers
 both $h$ and $g$. Hence $[u,v]$ is upper semimodular.
\end{proof}

\begin{proof}[Proof of Corollary~\ref{cor:Main}] 
By \cite[Theorem 3.1]{Bj}, a upper semilattice is shellable, hence Cohen-Macaulay.
\end{proof}

\subsection{Galois connection and proof of Theorem \ref{thm:mobius}}\label{sse:GaloisConnect}

For $I \subseteq J \subseteq S$, we will denote the poset
$$\PP_I(J) = \{ K \subseteq S | I \subseteq K \subseteq J \}$$
which is ordered by inclusion of subsets.

Choose $u,v \in W$ such that $u<v$ and set $I = \SS(u)$ and
$J=\SS(v)$.
We define two maps:
$$\GG_v : \PP_I(J) \longrightarrow [u,v]$$
where $\GG_v(K) = v_K$ and
$$\FF : [u,v] \longrightarrow \PP_I(J)$$
with $\FF(w) = \SS(w)$.

\begin{theorem} \label{th:maps} Let $u,v \in W$ such that $u \leq v$, and set $I = \SS(u)$ and
$J=\SS(v)$.
\begin{enumerate}
\item[$(a)$]  $\GG_v \circ \FF = id_{[u,v]}$ and $\FF\circ
\GG_v(K) \subseteq K$ for all $K \in \PP_I(J)$.
 \item[$(b)$] The map
$$\FF : ([u,v], \leq) \longrightarrow (\PP_I(J), \subseteq)$$
is a poset monomorphism. Moreover, for $w,g \in [u,v]$, $\FF(w
\vee g) = \FF(w) \cup \FF(g)$. 
\item[$(c)$] The map
$$\GG_v : (\PP_I(J), \subseteq, \cup, \cap) \longrightarrow ([u,v], \leq, \vee, \wedge)$$
is a poset epimorphism.  Moreover, $\GG_v( K \cap L ) = \GG_v(K) \wedge \GG_v(L)$
for $K, L \in \PP_I(J)$.
\item[$(d)$] For any $w \in [u,v]$,
$\GG_v^{-1}(w)$ is a sublattice of $(\PP_I(J), \subseteq, \cup,
\cap)$. Moreover,
$$\min \GG_v^{-1}(w) = \SS(w)\in \PP_I(J)$$ and
$$\max \GG_v^{-1}(w) = J \backslash \Des( v^{\SS(w)})\in \PP_I(J)$$
\end{enumerate}
\end{theorem}

Let $P = (P, \leq_P)$ and $Q = (Q, \leq_Q)$ be posets and $f : P
\rightarrow Q$ and $g : Q \rightarrow P$ be order preserving maps.
Recall that the pair $(f,g)$ is called a Galois connection if for
any $x \in Q$ and $y \in P$,
$$f(x) \leq_Q y \Leftrightarrow x \leq_P g(y).$$

\begin{cor} \label{cor:galois}
The pair $(\FF, \GG_v)$ is a Galois connection from $(\PP_I(J),
\subseteq)$ to $([u,v], \leq).$
\end{cor}
\begin{proof} 
By Theorem~\ref{th:maps} $\GG_v$ and $\FF$ are order preserving.
In addition, for $w \in [u,v]$ and $K \in \PP_I(J)$, if $w \leq
\GG_v(K)$, then by Theorem~\ref{th:maps} part $(d)$, $\FF(w) \subseteq
\FF(\GG_v(K)) \subseteq K$ and if $\FF(w) \subseteq K$, then $w =
\GG_v(\FF(w)) \leq  \GG_v(K)$.
\end{proof}

\begin{proof}[Proof of Theorem~\ref{th:maps}] (part $(a)$)
For a $w \in [u,v]$, we have that $\GG_v(\FF(w)) = \GG_v(\SS(w)) =
v_{\SS(w)} = w$ by definition of $w \leq v$.  For $K \in
\PP_I(J)$, then $\FF(\GG_v(K)) = \SS(v_K) \subseteq K$.

(part $(b)$)
$\FF$ is injective since $\GG_v \circ \FF = id_{[u,v]}$.  It is
order preserving by Lemma \ref{lemma:weakGasc}.  Lemma
\ref{lem:latticedef} part $(a)$ implies that $\FF( w \vee g ) =
\FF(w) \cup \FF(g)$.

 (part $(c)$)
$\GG_v$ is surjective since $\GG_v \circ \FF = id_{[u,v]}$. Take
$K, L$ such that $I \subseteq K \subseteq L \subseteq J$, then by
Lemma~\ref{lem:Useful}, $v_K \leq v_L$ and hence is order
preserving.


(part $(d)$)
Let $K, L \in \PP_I(J)$ be such that $\GG_v(K) = \GG_v(L) = w \in
[u,v]$. We have to show that $K \cap L$ and $K \cup L$ are in
$\GG_v^{-1}(w)$.  Note that $v_K = v_L = w$, therefore $w \in W_K
\cap W_L = W_{K \cap L}$.  Since $X_K, X_L \subseteq X_{K \cap
L}$, then $v_{K \cap L} = w$ because of the uniqueness of the
parabolic components. Hence $\GG_v( K \cap L) = w$.

From part $(a)$, $\SS(w) = \FF \circ \GG_v(K) \subseteq K$ which
implies that $\min \GG_v^{-1}(w) = \SS(w)$.
Moreover, since $v = v^K v_K = v^L v_L$ and $v_K = v_L = w$, we have that
$v^K = v^L \in X_K \cup X_L \subseteq X_{K \cup L}$.  Also $w \in
W_K \cup W_L \subseteq W_{K \cup L}$. Therefore $v_{K \cup L} = w$
and $K \cup L \in \GG_v^{-1}(w)$.
Let $A:=J \backslash \Des( v^{\SS(w)} )$. By definition  of $v^{\SS(w)}$,  we have
 $S(w)\cap \Des( v^{\SS(w)})=\emptyset$. Therefore   $\SS(w)\subseteq A$, 
and Lemma~\ref{lem:Useful} gives  $w \leq v_A$. 
Given that $I\subset \SS(w)$ we have $I\cap \Des( v^{\SS(w)})=\emptyset$ and 
 this gives that $A\subset \PP_I(J)$.  Since $v^{\SS(w)} \in  X_{A}$ and
$v_{\SS(w)} \in W_{\SS(w)} \subseteq W_A$, we have  $v^{\SS(w)} = v^A$ and
$w = v_A$, therefore $\GG_v(A) = w$.
\end{proof}

The value of the M\"obius function in Theorem~\ref{thm:mobius} follows from the
following classical result due to Rota.

\begin{prop} (\cite{R} Theorem 1)\label{prop:rota} If $(f,g)$ is a Galois connection between $P$ and $Q$, then for
$a \in P$ and $b \in Q$,
$$\sum_{\genfrac{}{}{0cm}{}{x \in P}{f(x) = b}}\mu_{P}(a,x) =
\sum_{\genfrac{}{}{0cm}{}{y \in Q}{g(y) = a}} \mu_{Q}(y,b) .$$
Moreover, both sums are equal to $0$ unless $g(f(a))= a$ and
$f(g(b))= b$, in which case they are both equal to
$\mu_{Q}(f(a),b) = \mu_{P}(a,g(b))$.
\end{prop}

\begin{proof}[Proof of Theorem \ref{thm:mobius}]
Corollary \ref{cor:galois} implies that Proposition
\ref{prop:rota} applies for $P = [u,v]$ and $Q = \PP_I(J)$ where
$b = J = S(v)$, $I = S(u)$, and $a=u$. Since $\FF : [u,v]
\rightarrow \PP_I(J)$ is injective, we have
\begin{equation}\label{eq:galmobius}
\mu(u,v)=\sum_{\genfrac{}{}{0cm}{}{w \in [u,v]}{\FF(w) =
J}}\mu_{[u,v]}(u,w) = \sum_{\genfrac{}{}{0cm}{}{L \in
\PP_I(J)}{\GG_v(L) = u}} \mu_{\PP_I(J)}(L,J) .
\end{equation}

Now Proposition \ref{prop:rota} says that $\mu(u,v) = 0$ unless
$\GG_v^{-1}(u)$ is a single element.  Theorem \ref{th:maps}
part $(d)$ identifies that $\GG_v^{-1}(u)$ will contain exactly
one element if and only if $S(u) = J \backslash \Des(u)$.  In this
case $\mu(u,v) = \mu(I,J) = (-1)^{|\SS(v)| - |\SS(u)|}$.

Now take $u=e$. Then the nonzero condition is equivalent to
 $\Des(v)=\SS(v)=J$, since $\Des(e)=\SS(e)=\emptyset$. This condition is a 
characterization of $w_0(J)$. 
\end{proof}

\begin{cor}
$$\sum_{x \leq v} \mu(e,x) t^{|\SS(x)|} = (1-t)^{|\Des(v)|}$$
\end{cor}
\begin{proof}
\begin{align*}
\sum_{x \leq v} \mu(e,x) t^{|\SS(x)|}
&=\sum_{\genfrac{}{}{0cm}{}{x \leq v}{x = w_0(\SS(x))}} (-1)^{|\SS(x)|} t^{|\SS(x)|}\\
&=\sum_{K \subseteq \Des(v)} (-1)^{|K|} t^{|K|}\\
&=(1-t)^{|\Des(v)|}
\end{align*}
\end{proof}

%
%
%
%

\section{Examples: Coxeter groups of type $A$, $B$ and $D$}\label{se:Examples}

\subsection{$W=\Sn_n$}

Let $\Sn_n^{(k)}$ denote the set of permutations $\{ w \in \Sn_n :
|\Gasc(w)| = k \}$. From the previous section, $(\Sn_n,
\leq)$ forms a graded poset with the elements of rank $k$ are
$\Sn_n^{(n-k)}$.  It is well known that $f_A(x) = 1-1/\sum_{n\geq
0} n! x^n = x + x^2 + 3 x^3 + 13 x^4 + 71 x^5 + \cdots$ is a
generating function for the number of elements in $\Sn_n^{(0)}$
(see \cite{C} or \cite[A003319]{OLEIS} and references therein). Now
any permutation can be identified with an ordered sequence of
permutations with empty connectivity set by breaking the
permutation at the positions of the elements in the connectivity
set and standardizing the segments. More precisely, if $\Gasc(u) =
\{ \tau_{i_1}, \tau_{i_2}, \ldots, \tau_{i_k} \}$, write $w=w_1\dots w_{k+1}$
 as in Equation (\ref{eq:ParabolicCompenent}), then
$$ \std_{\Gasc(u)}(w)=\std(w_1)\times\cdots\times\std(w_{k+1}) 
\in \Sn_{i_1}^{(0)} \times \Sn_{i_2-i_1}^{(0)} \times \cdots
\times \Sn_{n-i_k}^{(0)},$$ and this is clearly a bijection
between $\Sn_n$ and $\biguplus_\alpha \Sn_{\alpha_1}^{(0)} \times
\Sn_{\alpha_2}^{(0)} \times \cdots \times
\Sn_{\alpha_{\ell(\alpha)}}^{(0)}$ where the index $\alpha$ runs
over all compositions $(\alpha_1, \alpha_2, \ldots,
\alpha_{\ell(\alpha)})$ of $n$. This implies that the generating
function for $|\Sn_n^{(k)}|$ with $n\geq 1$ is $f_A(x)^{k+1}$ and
the proposition below follows from this remark (see also \cite{AS}
Corollary 6.4).

\begin{prop} \label{prop:typeAgf}
The coefficient of $x^n t^{k}$ in the generating function
$\frac{f_A(x)}{1-t f_A(x)}$ is equal to the number of permutations
in the set $\Sn_n^{(k)}$.
\end{prop}

\begin{align*}
\frac{f_A(x)}{1-t f_A(x)} &= 1 +  x + x^2 \left(t + 1\right) + x^3
\left(t^2 + 2\, t + 3\right)
 + x^4\left(t^3 + 3\, t^2 + 7\, t + 13\right) \\&\hskip .3in+
 x^5 \left(t^4 + 4\, t^3 + 12\, t^2 + 32\, t + 71\right) + \cdots
 \end{align*}

\subsection{$W = \B_n$}

Denote by $\Sn_{[\pm n]}$ the symmetric group acting on the set $(-[n])\cup[n]$.
 When $W$ is equal to the hyperoctahedral group $\B_n$ of order $2^n n!$
we can characterize the elements of this partial order more
combinatorially on on signed permutations.  Let $\B_n$ represent
the signed permutations: it is the subgroup of $\Sn_{[\pm n]}$
 consisisting of the element $w$ such that $w(-i)=-w(i)$ for all $i\in [n]$.
 As a Coxeter group, it is generated by the elements $\{ \tau_0, \tau_1, \ldots
\tau_{n-1} \}$ where $\tau_0$ is the transposition $(-1, 1)$ 
 and $\tau_i$ is the product of the transpositions $(i, i+1)$ and $(-i, -i-1)$.

\begin{prop}  For $u \in \B_n$ with $S=\{ \tau_0, \tau_1, \ldots
\tau_{n-1} \}$,
\begin{equation}
\GAsc(u) = \{ \tau_i : |u(j)| < u(k), \forall~0\leq j\leq i
< k< n \}.
\end{equation}
with the convention that $u(0) = 0$.
\end{prop}



\begin{proof}
Let $\tau_i \in \Gasc(u)$, then for $k > i$ and $j < i$,
$\tau_k \tau_j = \tau_j \tau_k$ and hence $u \in W_{\{\tau_0,
\tau_1, \ldots, \tau_{i-1}\}} \times W_{\{\tau_{i+1}, \ldots,
\tau_{n-1}\}} \cong \B_{i} \times \Sn_{n-i}$.  For this reason
$|u(j)| \leq i$ and $u(k) > i$ and hence $\tau_i \in \{ \tau_i :
|u(j)| < u(k), \forall 0 \leq j \leq i < k \leq n\}$.

Now assume that $\tau_i\in \{ \tau_i : |u(j)| < u(k), \forall 0
\leq j \leq i < k \leq n\}$, then for $k > i$, $u(k)>0$ and since
$u(k) > \max \{ |u(1)|, |u(2)|, \ldots, |u(i)|\}$, by the pigeon
hole principle we know that $u(k) > i$ (since it is larger than
$i$ different positive values).  For $j \leq i$, we know that
$|u(j)| < \min \{ u(i+1), \ldots, u(n) \} \leq i+1$ and hence $u \in
\B_i \times \Sn_{n-i}$ and hence $\tau_i \in \Gasc(u)$.
\end{proof}

Denote $\B_n^{(k)} =  \{ w \in \B_n : |\Gasc(w)| = k \}$. This
characterization of the connectivity set in type $B$  gives a
method for calculating the number elements of $\B_n^{(k)}$.

\begin{prop}
$$f_{B}(x) = \frac{\sum_{n \geq 0} 2^n n! x^n}{\sum_{n \geq 0} n! x^n}  = 1 + x + 5\, x^2 + 35\, x^3 + 309\, x^4 + 3287\, x^5 + 41005\, x^6
+ \cdots$$ is a generating function for the number of elements of
$\B_n$ with empty connectivity set. Moreover, the coefficient of
$x^n t^{k}$ in the generating function $f_B(x)/(1-t f_A(x))$ is
equal to the number of elements in the set $\B_n^{(k)}$.
\end{prop}

\begin{align*}
\frac{f_B(x)}{1-t f_A(x)} &= 1 + x\left(t + 1\right) + x^2
\left(t^2 + 2\, t + 5\right) + \
x^3 \left(t^3 + 3\, t^2 + 9\, t + 35\right)\\
&\hskip .2in + x^4 \left(t^4 + 4\, t^3 + 14\, t^2 + 56\, t +
309\right) +\cdots
\end{align*}
\begin{proof}
Let $r = \min \{ i : \tau_i \in \Gasc( u)\}$, then $u(1) u(2)
\cdots u(r)$ is a word of a signed permutation representing an
element with empty connectivity set and $std(u(r+1) u(r+2) \cdots
u(n))$ is a word of a permutation of size $n-r$.  It is not hard
to see that this operation defines a bijection between $\B_n$ and
the set $\biguplus_{r=0}^n \B_r^{(0)} \times \Sn_{n-r}$.
Therefore if $f_B(x) = \sum_{n \geq 0} |\B_n^{(0)}| x^n$, then the
generating function will satisfy $\sum_{n\geq 0} 2^n n! x^n =
f_B(x) \sum_{n\geq0} n! x^n$.
 Moreover,
the proof for Proposition \ref{prop:typeAgf} shows that  $f_B(x)
f_A(x)^k$ is a generating function for the number of elements of
$\B_n^{(k)}$.
\end{proof}

\begin{figure}[htbp] 
   \centering
   \includegraphics[width=3in]{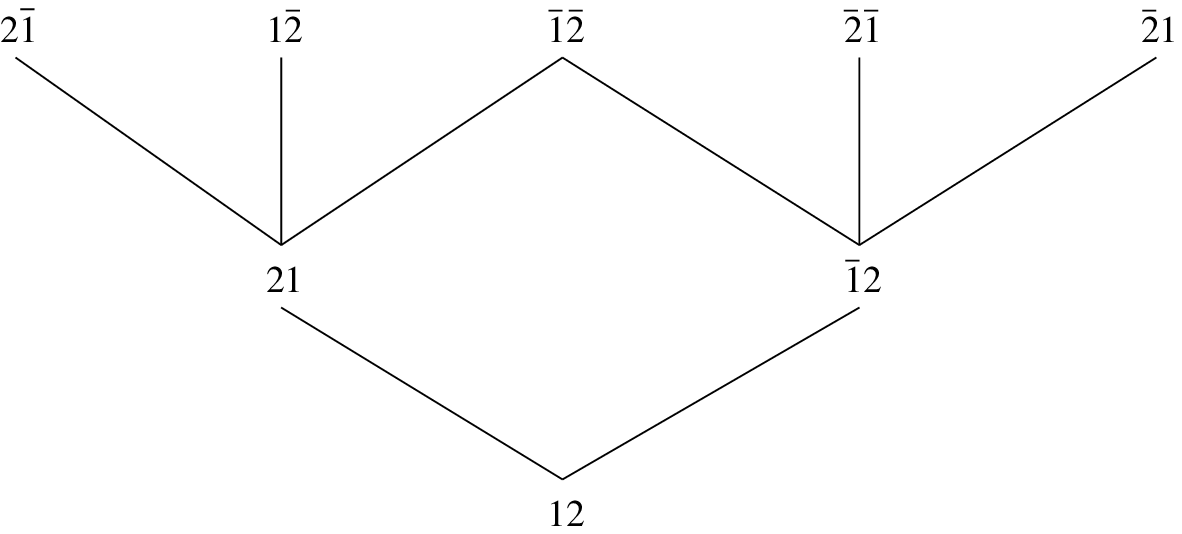}
   \caption{Hasse diagram for posets of $(\B_2, \leq)$}
   \label{hasseB2}
\end{figure}

\begin{figure}[htbp] 
   \centering
   \includegraphics[width=6in]{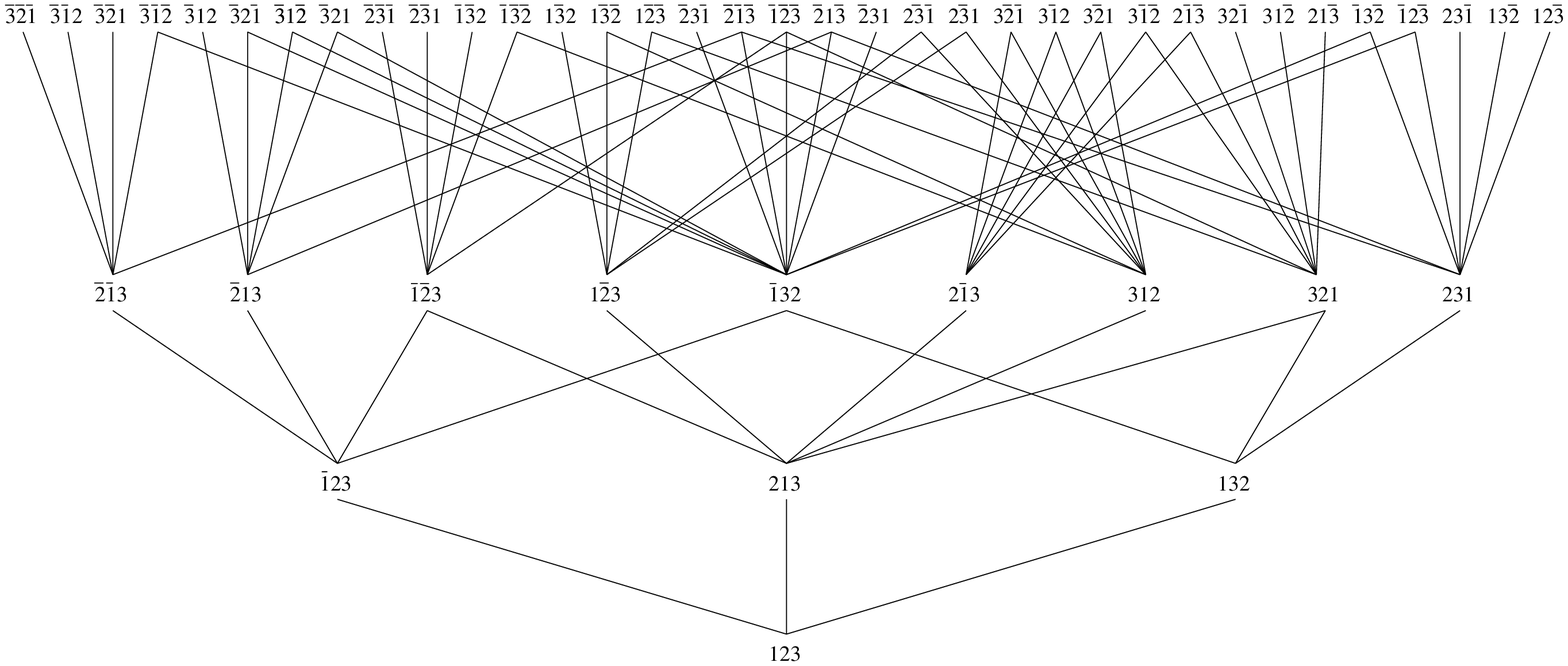}
   \caption{Hasse diagram for posets of $(\B_3, \leq)$}
   \label{hasseB3}
\end{figure}


\subsection{$W = \D_n$}

For $n\geq 2$, let $\D_n$ denote the Coxeter group of type $D$ of
order $2^{n-1} n!$. $\D_n$ can be realized as a subgroup of $\B_n$
of signed permutations and is generated by $\{ t_0, t_1, \ldots,
t_{n-1} \}$ where for $i\geq 1$, $t_i = \tau_i$ and $t_0= \tau_0 t_1 \tau_0$. In other words,
$t_0(1)=-2$ and $t_0(2)=-1$.

Denote $\D_n^{(k)} =  \{ w \in \D_n : |\Gasc(w)| = k \}$.

\begin{prop}  For $u \in \D_n$, if $(u(1) = 1$ and $u(k)>0$ for $1 \leq k \leq n )$ 
 or $( u(1) < -1$ and $u^{-1}(1) < -1$ and $|\{1\le i\le n: u(i)<0\}|=2 )$, then
\begin{equation}
\GAsc(u) = \{ t_1 \} \cup \{ t_i :  i \neq 1\hbox{ and }
|u(j)| < u(k), \forall~0\leq j\leq i < k< n\}
\end{equation}
otherwise
\begin{equation}
\GAsc(u) = \{ t_i : i \neq 1\hbox{ and } |u(j)| < u(k),
\forall~0\leq j\leq i < k< n \}.
\end{equation}
\end{prop}

\begin{proof} For $u\in\D_n\subset\B_n$, we denote $C_B(u)$ the connectivity
set of $u$ viewed as an element of $\B_n$.
Then observe that $\Gasc_{B}(u)\setminus\{t_1\}= \Gasc(u)\setminus\{t_1\}$.

We are left to consider when $t_1$ is in $\Gasc(u)$.
Clearly $\{t_0,t_1\}\subseteq C(w)$ if and only if $u(i)>0$ for all $i>0$ and $u(1) = 1$.
Consider now the case where $u \in W_{\{t_0, t_2, ..., t_{n-1}\}}$.
If $u \in W_{\{t_2...t_{n-1}\}}$, then $u(1) = 1$ and $u( k ) = 2$ for some  
$k\ge 2$. Hence
$t_0 u(1) = -2$ and $t_0 u(k) = -1$ and $|\{1\le i\le n: u(i)<0\}|=2$.
Assume by induction that
$u \in W_{\{t_0,t_2,t_3...t_{n-1}\}}$ and that the only two negative values of $u(i)$ are $u(1) = -\ell$
and $u(k) = -1$ for some $k\ge 2$.
If $\ell=2$, then $t_0 u$ has no negative signs and $t_0 u(1) = 1$ and hence
$t_0 u \in W_{\{t_2...t_{n-1}\}}$.
If $\ell\neq 2$, then $u(d) = 2$ for some $d$ and $t_0 u(d) = -1$
and $t_0 u(k) = 2$  and $t_0 u(1) = -\ell$.  Hence  by induction,
if $u \in W_{\{t_0, t_2, t_3, ... t_ {n-1}\}}$,
then $u(1) = 1$ and $u(k)>0$ for $1 \leq k \leq n$,
 or $u(1) < -1$ and $u^{-1}(1) < -1$ and $|\{1\le i\le n: u(i)<0\}|=2$.

For the converse, assume that $u(1) < -1$ and $u^{-1}(1) < -1$ and $|\{1\le i\le n: u(i)<0\}|=2$.
Then let $k = -u(1)$ and set $w = t_0 t_2 t_3 \cdots t_{k-1} u$.  It is easy to check
that $w \in \Sn_n$ and $w(1) = 1$.  Therefore $t_1 \in C(w)$ and so $t_1$ does not appear
in a reduced expression for $w$.
Therefore if $t_{i_1} \cdots t_{i_{\ell(w)}}$ is a reduced word for $w$, then
$t_{k-1} t_{k-2} \cdots t_2 t_0 t_{i_1} \cdots t_{i_{\ell(w)}}$ is a word for $u$
which does not contain a $t_1$ and so $t_1 \in C(u)$.
\end{proof}

\begin{prop}
\begin{equation}\label{eq:typeDgf}
f_D(x) = \frac{3 + \sum_{n\geq 0} 2^n n! x^n}{2 \sum_{n\geq 0} n!
x^n} + x -2 = x^2 + 13\, x^3 + 135\, x^4 + 1537\, x^5 + 19811\,
x^6 + \cdots
\end{equation}
is a generating function for the number of elements of $\D_n$ with
empty connectivity set. Moreover, the coefficient of $x^n t^{k}$
in the generating function $\frac{2 t f_A(x) + t^2 x f_A(x) +
f_D(x)}{(1-t f_A(x))}$ is equal to the number of elements in the
set $\D_n^{(k)}$.
\end{prop}

\begin{align*}
\frac{2 t f_A(x) -2tx + t^2 x f_A(x) + f_D(x)}{(1-t f_A(x))} =
~&x^2\, (t^2 + 2t + 1) + x^3\, \left(t^3 + 3\, t^2 +
 7\, t + 13\right)\\ &+ x^4\, (t^4 + 4\, t^3  + 12\, t^2 + 40\, t  + 135) + \cdots
\end{align*}
\begin{proof} Let $f_D(x) = \sum_{n \geq 2} | \D_n^{(0)} | x^n$ be the generating function
for the number of elements of $\D_n$ with empty connectivity set
and $f_A(x)$ be the generating function for the number of elements
of $\Sn_n$ with empty connectivity set from Proposition
\ref{prop:typeAgf}.

Now take $w$ to be an element of $\D_n^{(k)}$ such that neither
$t_0$ nor $t_1$ are elements in $\Gasc(w) = \{ t_{i_1},
\ldots, t_{i_k} \}$ then $i_1 > 1$ and $w \in W_{\{t_0, t_1,
\ldots, t_{i_1-1}\}} \times W_{\{t_{i_1+1}, \ldots, t_{i_2-1}\}}
\times \cdots \times W_{\{i_k+1, \ldots n-1\}} \cong \D_{i_1}
\times \Sn_{i_2-i_1} \times \cdots \times \Sn_{n-i_k}$.  Since
each component must have empty connectivity set, the generating
function for these elements is $f_D(x) f_A(x)^k$ .

Take $w$ to be an element of $\D_n^{(k)}$ such that both of $t_0$
or $t_1$ are elements in $\Gasc(w)$, then $w \in W_{\{t_2,
\ldots, t_{n-1}\}}$ and has exactly $k-2$ other elements in the
connectivity set. By Proposition  \ref{prop:typeAgf} the
generating function for these elements will be $x f_A(x)^{k-1}$.

Now for an element $w \in \D_n^{(k)}$ with exactly one of $t_0,
t_1 \in \Gasc(w)$ (take w.l.o.g. $t_0$), then $w \in W_{\{
t_1, t_2, \ldots, t_{n-1}\}} \cong \Sn_n$ and has exactly $k-1$
other elements in the connectivity set, but is not in $W_{\{ t_2,
\ldots, t_{n-1}\}}$. Therefore the generating function for these
elements is given by $f_A(x)^k - x f_A(x)^{k-1}$.

Since $\D_n = \biguplus_{k \geq 0} \D_n^{(k)}$, we have the
following generating function equation
$$\sum_{n \geq 2} 2^{n-1} n! x^n = 
f_D(x) + \sum_{k \geq 1} ( f_D(x) f_A(x)^k + x f_A(x)^{k-1} + 2( f_A(x)^k -
x f_A(x)^{k-1})) - x.$$ A bit of algebraic manipulation yields the
equation \eqref{eq:typeDgf} and
$$
\sum_{n \geq 2} |\D_n^{(k)}(x)| x^n t^k = f_D(x) + \sum_{k \geq 1}
t^k (f_D(x) f_A(x)^k + x f_A(x)^{k-1} + 2( f_A(x)^k - x
f_A(x)^{k-1})) - tx$$ is an expression for the bigraded generating
function formula in the proposition.
\end{proof}


\begin{thebibliography}{\hphantom{99}}

\bibitem{AS}   {\sc M.~Aguiar and F.~Sottile},
{\it Structure of the Malvenuto-Reutenauer Hopf
algebra of permutations}, Adv. Math {\bf 191} (2005), 225--275.

\bibitem{BZ} {\sc N.~Bergeron and M.~Zabrocki}, {\it The Hopf algebra of non-commutative symmetric functions
     and quasi-symmetric functions are free and cofree}, in preparation.

\bibitem{Bj}
{\sc A.~ Bj\"orner},  {\it Shellable and Cohen-Macaulay
partially ordered sets}, Trans. of the AMS {\bf 260(1)} (1980), 159--184.


\bibitem{bjorner-brenti} {\sc A.~Bj\"orner and F.~Brenti}, {\it Combinatorics of Coxeter Groups},
 Graduate Texts in Math. Springer (2005).

\bibitem{calan} {\sc D.~Callan}, {\it Counting stabilized-interval-free permutations},
 J.~Integer Sequences (electronic) {\bf 7} (2004), Article 04.1.8.

\bibitem{C} {\sc L.~Comtet}, {\em Advanced Combinatorics}, Reidel (1974).

\bibitem{NCSF6} {\sc G.~Duchamp, F.~Hivert and J.~Y.~Thibon} {\it Noncommutative symmetric functions VI: free quasi-symmetric functions and related algebras} , International Journal of Algebra and Computation {\bf 12} (2002), 671--717.


\bibitem{H} {\sc J.~E.~Humphreys}, {\it Reflection Groups and
Coxeter Groups}, Cambridge University Press, Cambridge, UK (1990).

\bibitem{R} {\sc G.-C.~Rota}, {\it On the foundations of combinatorial theory I: Theory of
M\"obius functions}, Z. Wahrscheinlichkeitstheorie {\bf 2} (1964), 340--368.  Reprinted
in ``Classic Papers in Combinatorics,'' (I. Gessel, G.C. Rota, Eds.), Birkh\"auser, Boston (1987).

\bibitem{St}   {\sc R.~Stanley}, {\it The Descent Set and Connectivity Set of a Permutation},
ArXive {\tt math.CO/0507224}.

\bibitem{St1}   {\sc R.~Stanley}, {\em Enumerative Combinatorics, Vol.~1},
Wadsworth and  Brooks/Cole, 1986.

\bibitem{St2} {\sc R.~Stanley}, {\em Enumerative Combinatorics, Vol. ~2},
Cambridge University Press, 1999.

\bibitem{matsumoto} {\sc H.~Matsumoto}, {\it G\'en\'erateurs et relations des groupes
 de Weyl g\'en\'eralis\'es}, C.~R.~Acad.~Sci.~Paris {\bf 258} (1964), 3419--3422.

\bibitem{OLEIS} {\sc N.~J.~A.~Sloane}, editor (2003),
The On-Line Encyclopedia of Integer Sequences, published electronically
at {\it http://www.research.att.com/\~{}njas/sequences/}.

\bibitem{tits} {\sc J.~Tits}, {\it Groupes et g\'eom\'etrie de Coxeter}, preprint I.H.E.S., Paris (1961).

\end{thebibliography}
\end{document}